\documentclass[12pt]{article}
\usepackage{amsfonts}
\usepackage{graphics}
\usepackage{graphicx}
\usepackage{amsmath,amssymb}

\usepackage{color}
\date{}

\usepackage[all]{xypic}

\begin{document}

\def\s{\subseteq}
\def\h{\widehat}
\def\v{\varphi}
\def\t{\widetilde}
\def\ov{\overline}
\def\L{\Lambda}
\def\O{\Omega}
\def\H{I\!\! H}
\def\a{\approx}
\def\la{\lambda}
\def\d{\delta}
\def\un{\underline}

\title{{\Large{\bf Admissible topologies on  $C(Y,Z)$ and ${\cal O}_Z(Y)$}}}

\author{Dimitris Georgiou$^a$, Athanasios Megaritis$^b$,\\ and Kyriakos Papadopoulos$^c$\\
\small{$^a$University of Patras, Department of Mathematics,  Greece}\\
\small{$^b$Technological Educational Institute of Messolonghi,}\\[-0.8ex]
\small{Department of Accounting, Greece}\\
\small{$^c$ University of Birmingham, School of Mathematics, United
Kingdom}}

\maketitle

\begin{abstract}
Let $Y$ and $Z$ be two given topological spaces, ${\cal O}(Y)$
(respectively, ${\cal O}(Z)$) the set of all open subsets of $Y$
(respectively, $Z$), and $C(Y,Z)$ the set of all continuous maps
from $Y$ to $Z$. We study Scott type topologies on ${\mathcal O}(Y)$
and we construct admissible topologies on $C(Y,Z)$ and ${\mathcal
O}_Z(Y)=\{f^{-1}(U)\in {\mathcal O}(Y): f\in C(Y,Z)\ {\rm and}\ U\in
{\mathcal O}(Z)\}$, introducing new problems in the field.
\end{abstract}


\medskip


\newtheorem{example}{Example}[section]
\newtheorem{remark}{Remark}[section]
\newtheorem{theorem}{Theorem}[section]
\newtheorem{corollary}[theorem]{Corollary}
\newtheorem{proposition}[theorem]{Proposition}
\newtheorem{definition}{Definition}
\newtheorem{question}[theorem]{Question}
\newtheorem{notation}[theorem]{Notation}

\section{Introduction and Preliminaries}

We denote by $Y$  and $Z$ two fixed topological spaces and by
$C(Y,Z)$ the set of all continuous maps from $Y$ to $Z$. If $t$ is a
topology on $C(Y,Z)$, then the corresponding topological space is
denoted by $C_{t}(Y,Z)$.

Let $X$ be a topological space and $Q$ a subset of $X$. By ${\rm
Cl}(Q)$ we denote the closure of $Q$ in $X$.

By ${\cal O}(X)$ we denote the family of all open subsets of  $X$
under a given topology and by ${\cal O}_Z(Y)$ the family
$\{f^{-1}(U):f\in C(Y,Z)\ \mbox{and}\  U\in {\cal O}(Z)\}$.

Let $X$ be a space and let $F:X\times Y\to Z$ be a continuous map.
By $\h{F}$ we denote the map from $X$ to the  set $C(Y,Z)$, such
that $\h{F}(x)(y)=F(x,y)$, for every $x\in X$ and $y\in Y$. Let $G$
be a map from $X$ to  $C(Y,Z)$. By $\t{G}$ we denote the map from
$X\times Y$ to  $Z$, such that $\t{G}(x,y)=G(x)(y)$, for every
$(x,y)\in X\times Y$.

\begin{definition}{\rm (see \cite{ARE}
and \cite{DUG}) A topology $t$ on $C(Y,Z)$ is called  {\it
admissible}, if for every space $X$, the continuity of the map
$G:X\to C_{t }(Y,Z)$ implies the continuity of the map
$\t{G}:X\times Y\to Z$ or equivalently the evaluation map
$e:C_{t}(Y,Z)\times Y\to Z$, defined by $e(f,y)=f(y)$ for every
$(f,y)\in C(Y,Z)\times Y$, is continuous.}
\end{definition}

\begin{definition}{\rm (see, for example,  \cite{SCO}) The {\it
Scott topology} $\Omega(Y)$ on ${\mathcal O}(Y)$ is defined as
follows: a subset $\H$ of ${\mathcal O}(Y)$ belongs to $\Omega(Y)$
if:

$(\alpha)$ $U\in \H$, $V\in {\mathcal O}(Y)$, and $U\s V$ imply
$V\in \H$ and

$(\beta)$ for every collection of open sets of $Y$, whose union
belongs to $\H$, there are finitely many elements of this collection
whose union  also belongs to $\H$.}
\end{definition}

\begin{definition}{\rm (see \cite{LAM})
The {\it strong Scott topology $\Omega_1 (Y)$} on ${\mathcal O}(Y)$
is defined as follows: a subset $\H$ of ${\mathcal O}(Y)$ belongs to
$\Omega _1(Y)$ if:

$(\alpha)$ $U\in \H$, $V\in {\mathcal O}(Y)$, and $U\s V$ imply
$V\in \H$ and

$(\beta)$ for every open cover  of $Y$, there are finitely many
elements of this cover whose  union  belongs to $\H$.}
\end{definition}

\begin{definition}{\rm (see, for example, \cite{LAM} and \cite{MN})
The {\it Isbell topology}
 on $C(Y,Z)$, denoted here by $t_{Is}$, is the topology which has as a subbasis the family
of all sets of the form:
$$(\H, U)=\{f\in C(Y,Z):f^{-1}(U)\in \H\},$$
\noindent where $\H\in\Omega(Y)$ and $U\in {\mathcal O}(Z)$.}
\end{definition}

\begin{definition}{\rm (see, for example, \cite{LAM} and \cite{MN})
The {\it strong Isbell topology}
 on $C(Y,Z)$, denoted here by $t_{sIs}$, is the topology which has as a subbasis the family
of all sets of the form:
$$(\H, U)=\{f\in C(Y,Z):f^{-1}(U)\in \H\},$$
\noindent where $\H\in\Omega_1(Y)$ and $U\in {\mathcal O}(Z)$.}
\end{definition}

\begin{definition}{\rm (see \cite{FOX})
The {\it compact open topology} on $C(Y,Z)$, denoted here by
$t_{co}$, is the topology  which has as a subbasis the family of all
sets of the form: $$(K,U)=\{f\in C(Y,Z):f(K)\subseteq U\},$$ where
$K$ is a compact subset of $Y$ and $U\in {\mathcal O}(Z)$.}
\end{definition}

It is known that $ t_{co}\subseteq t_{Is}$ (see, for example, \cite{MN}).

\begin{definition}{\rm (see, for example,
\cite{LAM}) A subset $K$ of  a space $X$ is said to be {\it
bounded}, if every open cover of $X$ has a finite subcover for $K$.}
\end{definition}

\begin{definition}{\rm (see, for example, \cite{SCO})
A space $X$ is called {\it corecompact}, if for every $x\in X$ and
for every open neighborhood $U$ of $x$ there exists an open
neighborhood $V$ of $x$ such that $V\subseteq U$ and the subset $V$
is bounded in the space $U$.}
\end{definition}

\begin{definition}{\rm (see, for example, \cite{LAM})
A space $X$ is called {\it locally bounded}, if for every $x\in X$
and for every open neighborhood $U$ of $x$ there exists a bounded
open  neighborhood $V$ of $x$ such that $V\subseteq U$.}
\end{definition}

Below, we  give some  well known results on admissible topologies:

(1) The  compact open topology on $C(Y,Z)$ is admissible, if $Y$ is
a regular locally compact space (see \cite{ARE}).

(2) The Isbell topology on $C(Y,Z)$ is admissible, if $Y$ is a
corecompact space
 (see, for example, \cite{LAM} and \cite{SW}).

(3) The strong Isbell topology on $C(Y,Z)$ is admissible, if $Y$ is
a locally bounded space
 (see, for example, \cite{LAM}).

(4) A topology which is larger than an admissible topology is also
admissible (see \cite{ARE}).

\medskip

For a summary of all the above results  and some open problems on function spaces see \cite{GEO}.

\medskip

In this paper, we give Scott type  topologies on the
set
${\mathcal O}(Y)$
 and
 define using these  topologies, in a standard way,
new admissible topologies on the sets $C(Y,Z)$ and ${\mathcal
O}_Z(Y)$. We finally introduce  questions on the field.

\section{Scott type topologies on the set ${\mathcal O}(Y)$}

Throughout the text, by $\tau_X$ we will denote the corresponding
topology on $X$, where $X$ is a topological space.

\begin{definition}
{\rm Let $Y$ and $Z$ be two topological spaces. The topology on $Y$,
denoted here by $\tau_Y^Z$, which has as a subbasis the family:
$${\mathcal O}_Z(Y)=\{f^{-1}(U)\in {\mathcal O}(Y): f\in C(Y,Z)\
{\rm and}\ U\in {\mathcal O}(Z)\}$$ is called the {\it $Z$-topology corresponding to the topology $\tau_Y$ of $Y$}.}
\end{definition}

Clearly, $\tau_Y^Z\s \tau_Y$.

\begin{example}
{\rm (1) Let $\mathbb R$ be the set of real numbers equipped with
its usual topology $\tau_{\mathbb R}$ and let $Z$ be any set
together with its trivial topology that is, the indiscrete topology. Then, ${\mathcal O}_Z(\mathbb
R)=\{\emptyset,\mathbb R\}$ and, therefore, $\tau_{\mathbb R}^Z\not
= \tau_{\mathbb R}$.

(2) It is well known that (see, for example, \cite{ST} (Example 92), \cite{IL}, and \cite{TZ})
for a fixed space $Z$ there exists a space $Y$ such that every continuous map from $Y$ to
$Z$ is constant. This means that ${\mathcal O}_Z(Y)=\{\emptyset,Y\}$ and, therefore,
$\tau_{Y}^Z\not = \tau_{Y}$.

(3) Let $Z={\bf S}$ be the Sierpi\'{n}ski space, that is  $Z=\{0,1\}$
and $\tau_Z=\{\emptyset,\{1\},\{0,1\}\}$. If $Y$ is another
topological space, then $C(Y,{\bf S})=\{ {\mathcal X}_V : V \in
{\mathcal O}(Y)\}$, where ${\mathcal X}_V : Y \to {\bf S}$ denotes the
characteristic function of $V$,
$${\mathcal X}_V(y)=\begin{cases}1 & {\rm if } \ y\in V,\\
                          0 & {\rm if } \ y\not\in V.\end{cases}$$
                          In this case we observe that:
\begin{eqnarray*}
 {\mathcal O}_Z(Y) &=& \{{\mathcal X}_V^{-1}(\{1\}): {\mathcal X}_V \in C(Y,{\bf S})\}\\
                   &=& \{V : V\in {\mathcal O}(Y) \}\\
                   &=& {\mathcal O}(Y).
                   \end{eqnarray*}
  Thus, $\tau_Y^{\bf S} = \tau_Y$.}
\end{example}

\begin{definition}
{\rm Let $Y$ and $Z$ be two topological spaces. A subset $K$ of $Y$
is called {\it $Z$-compact}, if $K$ is compact in the space
$(Y,\tau_Y^Z)$, where $\tau_Y^Z$ is the $Z$-topology corresponding to the topology $\tau_Y$ of $Y$.}
\end{definition}

\begin{example}
{\rm  Every compact subset of a space  $Y$ is $Z$-compact, but the
converse statement is not true. Indeed, let $\mathbb R$ be the set
of real numbers with the usual topology and $Z$ a set equipped with
the trivial topology. Then, every subset of $\mathbb R$ is
$Z$-compact, while subsets of the set of real numbers are not
necessarily compact, in general.}
\end{example}

\begin{definition}
{\rm Let $Y$ and $Z$ be two topological spaces. The {\it $Z$-compact
open topology} on $C(Y,Z)$, denoted here by $t_{co}^Z$, is the
topology, which has as a subbasis the family of all sets of the
form:
$$(K,U)=\{f\in C(Y,Z):f(K)\subseteq U\}, $$
where $K$ is a $Z$-compact subset of $Y$ and $U\in {\mathcal O}(Z)$.}
\end{definition}

\begin{remark}
{\rm For the  topologies   $t_{co}$ and $t_{co}^Z$ we have that:
 $$t_{co}\s t_{co}^Z.$$}
 \end{remark}

\begin{example}
{\rm Let  $Y$ be an arbitrary topological space and  ${\bf S}$ the Sierpi\'{n}ski space. By Example 2.1(3), we have that a subset $K$ of
$Y$ is ${\bf S}$-compact if and only if $K$ is compact in
$(Y,\tau_Y)$. In this case, we have  $t_{co}^{\bf S}=t_{co}$.

Furthermore, we observe that,  whenever we consider  the topology on
${\mathcal O}(Y)$ which has as subbasis the family:
$$\{<K>: K \textrm{ is compact in } Y\},$$ where  $<K> = \{U \in {\mathcal O}(Y) : K \subset U\}$,
then the topological spaces $C_{t_{co}^{\bf S}}(Y,{\bf S})$ and
${\mathcal O}(Y)$ are homeomorphic.  Indeed, it would be enough to
consider the homeomorphism $T:C(Y,{\bf S})\to {\mathcal O}(Y)$ for
which $T({\mathcal X}_V)=V$, for every  ${\mathcal X}_V \in C(Y,{\bf
S})$.

So, in the case where $Z={\bf S}$, it would be that $t_{co}^{\bf
S}=t_{co}$ on $C(Y,{\bf S})$ and this topology coincides with the topology
on ${\mathcal O}(Y)$, which has as subbasis the set $\{<K>: K
\textrm{ is compact in } Y\}$.}
\end{example}

\begin{definition}
{\rm By $\tau_1^Z$  we denote the family of all subsets of
${\mathcal O}(Y)$ that are defined as follows: a subset $\H $ of
${\mathcal O}(Y)$ belongs to $\tau_1^Z$ if:

$(\alpha)$ $f^{-1}(U)\in \H\cap {\mathcal O}_Z(Y)$, $V\in {\mathcal
O}(Y)$ and $f^{-1}(U)\s V$ imply $V\in \H $ and

$(\beta)$ for every collection
$\{f_{\lambda}^{-1}(U_{\lambda}):\lambda\in \Lambda\}$ of elements
of ${\mathcal O}_Z(Y)$, whose union belongs to $\H $, there are
finitely many elements $f^{-1}_{\lambda_i}(U_{\lambda_i}),$
$i=1,2,\cdots,n$ of this collection, such that:
$$\bigcup\{f^{-1}_{\lambda_i}(U_{\lambda_i}):i=1,2,\cdots,n\}\in\H.$$}
\end{definition}

\begin{proposition}
{\rm The family $\tau_1^Z$ defines a topology on ${\mathcal O}(Y)$,
called the {\it $Z$-Scott topology}.}
\end{proposition}

{\sl Proof.} The proof  follows trivially by Definition 13. $\Box$

\begin{definition}
{\rm By   $\tau_{1,s}^Z$ we denote the family of all subsets of
${\mathcal O}(Y)$ that are defined as follows: a subset $\H $ of
${\mathcal O}(Y)$ belongs to $\tau_{1,s}^Z$ if:

$(\alpha)$ $f^{-1}(U)\in \H\cap {\mathcal O}_Z(Y)$, $V\in {\mathcal
O}(Y)$ and $f^{-1}(U)\s V$ imply $V\in \H $ and

$(\beta)$ for every collection
$\{f_{\lambda}^{-1}(U_{\lambda}):\lambda\in \Lambda\}$ of elements
of ${\mathcal O}_Z(Y)$ whose union is equal to $Y$, there are
finitely many elements $f^{-1}_{\lambda_i}(U_{\lambda_i}),$
$i=1,2,\cdots,n$ of this collection, such that:
$$\bigcup\{f^{-1}_{\lambda_i}(U_{\lambda_i}):i=1,2,\cdots,n\}\in\H.$$}
\end{definition}

\begin{proposition}
{\rm The family $\tau_{1,s}^Z$ defines a topology on ${\mathcal
O}(Y)$, called the {\it strong $Z$-Scott topology}.}
\end{proposition}

{\sl Proof.} The proof  follows trivially by Definition 14. $\Box$

\begin{example}
{\rm (1) Let $Y=\{0,1\}$ be  equipped with the topology
$\tau_{Y}=\{\emptyset,\{0\},Y\}$ and let $Z$ be a set equipped with
the trivial topology. Then, $${\mathcal O}_Z(Y)=\{\emptyset,Y\},$$
$$\Omega(Y)=\Omega_1(Y)=\{\emptyset,\{Y\},\{\emptyset,Y\},\{\{0\},Y\},\{\emptyset,\{0\},Y\}\}$$
and
$$\tau_1^Z=\tau_{1,s}^Z=\{\emptyset,\{\{0\}\},\{Y\},\{\emptyset,Y\},\{\{0\},Y\},\{\emptyset,\{0\},Y\}\}.$$

(2) If $Y$ is an arbitrary topological space and ${\bf S}$ is the Sierpi\'{n}ski space, then a subset  $\mathbb{H}$ of ${\mathcal O}(Y)$,
will belong to $\tau_1^{\bf S}$, if:

\begin{enumerate}
\item ${\mathcal X}_V^{-1}(\{1\}) \in \mathbb{H} \cap {\mathcal O}_{\bf S}(Y)$, $V \in {\mathcal O}(Y)$, and ${\mathcal X}_V^{-1}(\{1\}) \subset W$ implies that
$W \in \mathbb{H}$. Equivalently, $V \in \mathbb{H} \cap {\mathcal
O}(Y) = \mathbb{H}$, $V \in {\mathcal O}(Y)$ and $V \subset W$
implies that $W \in H$.

\item For every collection $\{{\mathcal X}_V^{-1}(\{1\}): V \in {\mathcal O}(Y)\}$ of elements of ${\mathcal O}_{\bf S}(Y) = {\mathcal O}(Y)$, whose union belongs
to $\mathbb{H}$, there exist finitely many elements ${\mathcal
X}_{V_i}^{-1}(\{1\})$, $i=1,2,\cdots,n$, such that:
$$\bigcup\{{\mathcal X}_{V_i}^{-1}(\{1\}) : i = 1,2,\cdots,n\} \in
\mathbb{H}.$$ Equivalently, for every collection $\{V_{\lambda}\in
{\mathcal O}(Y):\lambda\in\Lambda\}$, such that $\bigcup\{V_{\lambda}:\lambda\in\Lambda\}\in \mathbb{H}$, there are
finitely many elements $\{V_{\lambda_i}:i =1,2,\cdots,n\}$, such that
$\bigcup\{V_{\lambda_i}:i =1,2,\cdots,n\}\in \mathbb{H}$. So, $\tau_1^Z$ is the Scott
topology on ${\mathcal O}(Y)$.
\end{enumerate}}
\end{example}

\begin{remark}
{\rm For the topologies $\Omega (Y)$, $\Omega_1(Y)$, $\tau_1^Z $,
and  $\tau_{1,s}^Z$ we have the following comparison:
$$\begin{array}{ccc}
\vspace{0.2cm}
\Omega_1 (Y) & \subseteq & \tau_{1,s}^Z\\
\vspace{0.2cm}
\reflectbox{\rotatebox[origin=c]{-90}{$\supseteq$}} &  & \reflectbox{\rotatebox[origin=c]{-90}{$\supseteq$}}\\
\Omega (Y)  & \subseteq   & \tau^Z_{1}\end{array}$$
}
\end{remark}

\begin{definition}
{\rm The $t_1^Z$ topology on $C(Y,Z)$ is the topology which has as a
subbasis the family of all sets of the form:
$$(\H, U)=\{f\in C(Y,Z):f^{-1}(U)\in \H\},$$
\noindent
where $\H$ is open in the  topology  $\tau_1^Z$ on ${\mathcal O}(Y)$ and $U\in {\mathcal O}(Z)$.}
\end{definition}

\begin{definition}
{\rm The $t_{1,s}^Z$ topology  on $C(Y,Z)$ is the topology which has
as a subbasis the family of all sets of the form:
$$(\H, U)=\{f\in C(Y,Z):f^{-1}(U)\in \H\},$$
\noindent
where $\H$ is open in the  topology  $\tau_{1,s}^Z$ on ${\mathcal O}(Y)$
 and $U\in {\mathcal O}(Z)$.}
\end{definition}

\begin{proposition}
{\rm Let $Y$ and $Z$ be two topological spaces and let $K$ be a
$Z$-compact subset of $Y$. Then, the set:
$$\H_K=\{U\in {\mathcal O}(Y):K\subseteq U\}$$ is open in ${\mathcal O}(Y)$ with the topology $\tau_1^Z$.}
\end{proposition}

{\sl Proof.} Let $f^{-1}(U)\in \H_K\cap {\mathcal O}_Z(Y)$, $V\in
{\mathcal O}(Y)$ and let $f^{-1}(U)\subseteq V$. Then, $K\subseteq
f^{-1}(U)\subseteq V$ and, therefore, $V\in \H_K$.

Now, let $\{f_{\lambda}^{-1}(U_{\lambda}):\lambda\in \Lambda\}$ be a
collection of sets of ${\mathcal O}_Z(Y)$, whose union belongs to
$\H_K$. Then:
$$K\subseteq \bigcup\{f_{\lambda}^{-1}(U_{\lambda}):\lambda\in \Lambda\}.$$
Since $K$ is $Z$-compact,  there are finitely many elements
$f^{-1}_{\lambda_i}(U_{\lambda_i}),$ $i=1,2,\cdots,n$ of this
collection such that:
$$K\subseteq \bigcup\{f^{-1}_{\lambda_i}(U_{\lambda_i}):i=1,2,\cdots,n\}$$
and, therefore,
$$\bigcup\{f^{-1}_{\lambda_i}(U_{\lambda_i}):i=1,2,\cdots,n\}\in \H_K.$$
Thus, the set $\H_K$ is open in ${\mathcal O}(Y)$ with the topology $\tau_1^Z$. $\Box$

\begin{remark}
{\rm  (1)  We observe that for every $Z$-compact subset $K$ of $Y$
we have:
\begin{eqnarray*}
(\H_K,U)&=&\{f\in C(Y,Z):f^{-1}(U)\in \H_K\}\\
&=&\{f\in C(Y,Z):K\subseteq f^{-1}(U)\}\\
&=&\{f\in C(Y,Z): f(K)\subseteq U\}.
\end{eqnarray*}
This says that $t^Z_{co}\subseteq t^Z_1$. Thus, by
 Remarks 2.1 and 2.2   we get the following comparison between the
topologies   $t_{co}$, $t_{co}^Z$, $t_{Is}$, $t_{sIs}$,
 $t_1^Z$, and $t_{1,s}^Z$:
$$\begin{array}{ccccc}
\vspace{0.2cm}
t_{co}^Z   & \subseteq & t_{1}^Z  & \subseteq & t_{1,s}^Z\\
\vspace{0.2cm}
\reflectbox{\rotatebox[origin=c]{-90}{$\supseteq$}}  &  & \reflectbox{\rotatebox[origin=c]{-90}{$\supseteq$}}  & & \reflectbox{\rotatebox[origin=c]{-90}{$\supseteq$}}  \\
t_{co}  & \subseteq   & t_{Is}  & \subseteq & t_{sIs}
\end{array}$$

(2) Let $${\mathcal T} = \{\tau_Y^Z : Z\ \mbox{is an arbitrary
topological space}\}.$$ We immediately see that $({\mathcal
T},\subseteq)$ has an upper bound, namely $\tau_Y$, which is also
the maximal element for ${\mathcal T}$ because, if $Z={\bf S}$, then
$\tau_Y = \tau_Y^{\bf S}$.

 In a similar way, we can prove that the set:  $${\mathcal T}_{co}=\{t_{co}^Z: Z\ \textrm{ is an arbitrary topological space}\}$$
 has a lower bound, namely $t_{co}$, which is a minimal element and, if
$Z={\bf S}$, then $t_{co}=t_{co}^{\bf S}$.

Also, in a similar manner  the set:
$${\mathcal T}_{1}=\{t_1^Z: Z\ \textrm{ is an arbitrary topological space}\}$$ has a
lower bound, namely $t_{Is}$, which is a minimal element and, if
$Z={\bf S}$, then $t_{Is}=t_{1}^{\bf S}$.

Finally, the set:  $${\mathcal T}_{1,s}=\{t_{1,s}^Z: Z\ \mbox{is an
arbitrary topological space}\}$$ has a lower bound, namely
$t_{sIs}$, which is a minimal element and, if $Z={\bf S}$, then
$t_{sIs}=t_{1,s}^{\bf S}$.}
\end{remark}

\begin{theorem}
{\rm Let $Z$ be a ${\rm T}_i$-space, where $i=0,1,2$. Then, the topological spaces $C_{t_{co}^Z}(Y,Z)$, $C_{t_{1}^Z}(Y,Z)$, and $C_{t_{1,s}^Z}(Y,Z)$ are also ${\rm T}_i$-spaces.}
\end{theorem}

{\sl Proof.} Since $Z$ is a ${\rm T}_i$-space, where $i=0,1,2$, the
space $C_{t_{co}}(Y,Z)$ will also  be a ${\rm T}_i$-space (see, for
example, \cite{DUG}). Thus, by Remark 2.3, the  spaces
$C_{t_{co}^Z}(Y,Z)$, $C_{t_{1}^Z}(Y,Z)$ and $C_{t_{1,s}^Z}(Y,Z)$ are
${\rm T}_i$-spaces. $\Box$

\begin{theorem}
{\rm The following statements are true:

(1) If $Y$ is a regular locally compact space, then the topologies
$t^Z_{co}$, $t^Z_1$ and $t^Z_{1,s}$ are admissible.

(2) If $Y$ is a corecompact space, then the topologies
 $t^Z_1$ and $t^Z_{1,s}$ are admissible.

(3) If $Y$ is a locally bounded  space, then the topology
 $t^Z_{1,s}$ is admissible.}
 \end{theorem}

 {\sl Proof.} The proof of this theorem follows from Remark 2.3
 and from the fact that a topology larger than an admissible topology is also
 admissible. $\Box$

\begin{definition}
{\rm Let $Y$ and $Z$ be two topological spaces. The space $Y$ is
called {\it locally $Z$-compact}, if the space $(Y,\tau_Y^Z)$, where
$\tau_Y^Z$ is the $Z$-topology corresponding to the topology
$\tau_Y$ of $Y$, is locally compact.}
\end{definition}

\begin{remark}
{\rm We observe that, if a space $Y$ is regular locally $Z$-compact,
then for every $y\in Y$ and for every neighborhood $f^{-1}(W)\in \tau_Y^Z$ of $y$, where
$f\in C(Y,Z)$ and $W\in {\mathcal O}(Z)$, there exists $g^{-1}(V)\in
\tau_Y^Z$ such that the closure of the set $g^{-1}(V)$  is compact
in the space $(Y,\tau_Y^Z)$ and $y\in g^{-1}(V)\subseteq {\rm
Cl}(g^{-1}(V))\subseteq f^{-1}(W)$.}
\end{remark}

\begin{theorem}
{\rm Let $Z$ be a space and $Y$ a regular locally $Z$-compact  space. Then,
 the $t^Z_{co}$ topology  on $C(Y,Z)$ is admissible.}
\end{theorem}

{\sl Proof.} It is sufficient to prove that the evaluation
map:
$$e:C_{t^Z_{co}}(Y,Z)\times Y\to Z$$ is continuous.

Let $(f,y)\in C(Y,Z)\times Y$ and let also $W\in {\mathcal O}(Z)$,
such that:
$$e(f,y)=f(y)\in W.$$ Then, we have that:
$$y\in f^{-1}(W).$$ Since the space $Y$ is regular locally $Z$-compact, there exists $g^{-1}(V)\in
\tau_Y^Z$, such that the set ${\rm Cl}(g^{-1}(V))$ is compact in the
space $(Y,\tau_Y^Z)$ and
$$y\in g^{-1}(V)\s {\rm Cl}(g^{-1}(V)) \s f^{-1}(W).$$ Since ${\rm Cl}(g^{-1}(V)) \s f^{-1}(W)$, we have that $f\in ({\rm Cl}(g^{-1}(V)),W
)$. Thus, $$({\rm Cl}(g^{-1}(V)), W)\times g^{-1}(V)$$ is an open
neighborhood of $(f,y)$ in $C_{t_{co}^Z}(Y,Z)\times Y$.

We finally prove that:
$$e(({\rm Cl}(g^{-1}(V)), W)\times g^{-1}(V))\s W.$$ Let $(h,z)\in ({\rm Cl}(g^{-1}(V)), W)\times g^{-1}(V)$. Then:
$$h\in ({\rm Cl}(g^{-1}(V)), W) \ {\rm and} \ z\in g^{-1}(V).$$
Therefore, $h(g^{-1}(V))\subseteq h({\rm Cl}(g^{-1}(V)))\subseteq W$ and $e(h,z)=h(z)\in W$.

Thus, the evaluation map $e$ is continuous and, therefore, the
$t_{co}^Z$ topology on $C(Y,Z)$ is admissible too. $\Box$

\begin{theorem}
{\rm Let $X$, $Y$, and $Z$ be three topological spaces. If the space
$Y$ is locally $Z$-compact, then the map: $$T: C_{t^Y_{co}}(X,Y)
\times C_{t^Z_{co}}(Y,Z) \to C_{t^Z_{co}}(X,Z),$$ with
$T(f,g)=g\circ f$ for every $(f,g) \in C(X,Y)\times C(Y,Z)$, is
continuous.}
\end{theorem}

{\sl Proof.} Let $(f,g) \in C(X,Y) \times C(Y,Z)$ and let: $$T(f,g)
= g \circ f \in (K,U)\in t^Z_{co},$$
 where $K$ is a compact subset of the space $(X,\tau^Z_X)$ and  $U$ an open subset of $Z$.
Then, we have that:
$$(g \circ f)(K)\subseteq U$$
or, equivalently:
$$K\subseteq f^{-1}(g^{-1}(U)).$$
Now, since the space $Y$ is locally $Z$-compact, the space
$(Y,\tau^Z_Y)$ is locally compact. Thus, for an arbitrary $y\in
g^{-1}(U)$ there exists $W_y\in {\mathcal O}(Z)$  such that ${\rm
Cl}(W_y)$ is compact in  $(X,\tau^Z_X)$ and
$$y\in W_y\subseteq {\rm Cl}(W_y) \subseteq g^{-1}(U). \eqno{(1)}$$
So, we have  that: $$g^{-1}(U)=\bigcup\{W_y:y\in g^{-1}(U)\}$$ and,
therefore,
$$K\subseteq f^{-1}(g^{-1}(U))=\bigcup\{f^{-1}(W_y):y\in g^{-1}(U)\}.$$
Since $K$ is a compact subset of the space $(X,\tau^Z_X)$, there are
finitely many elements $y_1,\cdots,y_n\in g^{-1}(U)$, such that:
$$K\subseteq \bigcup\{f^{-1}(W_{y_i}):i\in\{1,\cdots,n\}\}.$$ So, we have
$$f(K)\subseteq \bigcup\{W_{y_i}:i\in\{1,\cdots,n\}\}$$ or, equivalently: $$f\in (K,
\bigcup\{W_{y_i}:i\in\{1,\cdots,n\}\}).$$ Also, by relation (1) we have
$$g(\bigcup\{{\rm Cl}(W_{y_i}):i\in\{1,\cdots,n\}\})\subseteq U$$ and,
therefore, $$g\in(\bigcup\{{\rm Cl}(W_{y_i}):i\in\{1,\cdots,n\}\},U).$$
We observe that: $$(K, \bigcup\{W_{y_i}:i\in\{1,\cdots,n\}\})\in
t^Y_{co}$$ and that
$$(\bigcup\{{\rm Cl}(W_{y_i}):i\in\{1,\cdots,n\}\},U)\in t^Z_{co}.$$ By
all the above it suffices to prove that: $$T\Big((K,
\bigcup\{W_{y_i}:i\in\{1,\cdots,n\}\})\times(\bigcup\{{\rm
Cl}(W_{y_i}):i\in\{1,\cdots,n\}\},U)\Big)\subseteq (K,U).$$ Let
$$(h_1,h_2)\in (K, \bigcup\{W_{y_i}:i\in\{1,\cdots,n\}\})\times(\bigcup\{{\rm
Cl}(W_{y_i}):i\in\{1,\cdots,n\}\},U).$$ Then, $$h_1(K)\subseteq
\bigcup\{W_{y_i}:i\in\{1,\cdots,n\}\}$$ and $$h_2(\bigcup\{{\rm
Cl}(W_{y_i}):i\in\{1,\cdots,n\}\})\subseteq U.$$ Therefore,
\begin{eqnarray*}
(h_2\circ h_1)(K)=h_2(h_1(K))&\subseteq&
h_2(\bigcup\{W_{y_i}:i\in\{1,\cdots,n\}\})\\
&\subseteq& h_2(\bigcup\{{\rm
Cl}(W_{y_i}):i\in\{1,\cdots,n\}\})\subseteq U,
\end{eqnarray*}
so that $T(h_1,h_2)=h_2\circ h_1\in (K,U)$. Thus, the map $T$ is
continuous. $\Box$

\begin{definition}
{\rm Let $Y$ and $Z$ be two topological spaces. A subset $B$ of $Y$
is called {\it $Z$-bounded}, if $B$ is bounded in the space
$(Y,\tau_Y^Z)$, where $\tau_Y^Z$ is the $Z$-topology corresponding to the topology $\tau_Y$ of $Y$.}
\end{definition}

\begin{definition}
{\rm Let $Z$ be a space. A space $Y$ is called {\it locally
$Z$-bounded}, if for every $y\in Y$ and for every open neighborhood
$U$ of $y$, there exists a $Z$-bounded neighborhood $g^{-1}(V)\in
{\mathcal O}_Z(Y)$ of $y$, such that $g^{-1}(V)\s U$.}
\end{definition}

\begin{theorem}
{\rm Let $Z$ be a space and let $Y$ be a locally $Z$-bounded space.
Then, the $t^Z_{1,s}$ topology on $C(Y,Z)$ is admissible.}
\end{theorem}

{\sl Proof.} It is sufficient to prove that the evaluation map:
$$e:C_{t^Z_{1,s}}(Y,Z)\times Y\to Z$$ is continuous.

For this, we let $(f,y)\in C(Y,Z)\times Y$ and let $W\in {\mathcal
O}(Z)$, such that:
$$e(f,y)=f(y)\in W.$$ Then:
$$y\in f^{-1}(W).$$ Also, since  $Y$ is locally $Z$-bounded, there exists a
$Z$-bounded neighborhood  $g^{-1}(V)\in {\mathcal O}_Z(Y)$ of $y$,
such that:
$$y\in g^{-1}(V)\s f^{-1}(W).$$ Consider the set:
$$\H_{g^{-1}(V)}=\{U\in {\mathcal O}(Y): g^{-1}(V)\s U\}.$$
We prove that $\H_{g^{-1}(V)}$ belongs to the $\tau^Z_{1,s}$
topology.

Indeed, let $h^{-1}(U)\in \H_{g^{-1}(V)}\cap {\mathcal O}_Z(Y)$,
$U_1\in {\mathcal O}(Y)$ and $h^{-1}(U)\s U_1$. Then:
$$g^{-1}(V)\s h^{-1}(U)
\s U_1$$ and, therefore, $U_1\in \H_{g^{-1}(V)}$. Now, let
$\{f^{-1}_i(U_i):i\in I\}$ be a collection of elements of ${\mathcal
O}_Z(Y)$, whose union is equal to the set $Y$. Then, for every $y\in
Y$, there exists $i_y\in I$, such that $y\in f^{-1}_{i_y}(U_{i_y})$.
Since $Y$ is locally $Z$-bounded,  there exists a $Z$-bounded
neighborhood $g^{-1}_{i_y}(V_{i_y})\in {\mathcal O}_Z(Y)$ of $y$,
such that:
$$y\in g^{-1}_{i_y}(V_{i_y})\s f^{-1}_{i_y}(U_{i_y}).$$ We can easily deduce that the set
$\{g^{-1}_{i_y}(V_{i_y}):i_y\in I\}$ is an open cover of
$(Y,\tau_Y^Z)$. Since $g^{-1}(V)$ is $Z$-bounded, there exist
finitely many sets $$g^{-1}_{i_{y_1}}(V_{i_{y_1}}),\cdots,
g^{-1}_{i_{y_n}}(V_{i_{y_n}})$$ such that:
$$g^{-1}(V)\s \bigcup\{g^{-1}_{{y_{i_k}}}(V_{{y_{i_k}}}):k=1,2,\cdots,n\}.\eqno(2)$$ We now consider the sets
$f^{-1}_{i_{y_1}}(U_{i_{y_1}}),\cdots, f^{-1}_{i_{y_n}}(U_{i_{y_n}})$
of ${\mathcal O}_Z(Y)$, for which:
$$g^{-1}_{i_{y_k}}(V_{i_{y_k}})\s f^{-1}_{i_{y_k}}(U_{i_{y_k}}), \ \ k=1,2,\cdots,n.$$ But relation (2)
gives:
$$ g^{-1}(V)\s \bigcup\{g^{-1}_{i_{y_k}}(V_{i_{y_k}}):k=1,2,\cdots,n\}\s\bigcup\{f^{-1}_{i_{y_k}}(U_{i_{y_k}}):k=1,2,\cdots,n\}$$
and, therefore:
$$\bigcup\{f^{-1}_{i_{y_k}}(U_{i_{y_k}}):k=1,2,\cdots,n\}\in
\H_{g^{-1}(V)}.$$ So, the set $\H_{g^{-1}(V)}$ is open in the
$\tau_{1,s}^Z$ topology.

Also, since $g^{-1}(V)\s f^{-1}(W)$, we have that $f\in
(\H_{g^{-1}(V)},W )$. Thus, the set $(H_{g^{-1}(V)}, W)\times
g^{-1}(V)$ is an open neighborhood of $(f,y)$ in
$C_{t_{1,s}^Z}(Y,Z)\times Y$.

We finally prove that:
$$e((\H_{g^{-1}(V)}, W)\times g^{-1}(V))\s W.$$ Let $(h,z)\in (\H_{g^{-1}(V)}, W
)\times g^{-1}(V)$. Then:
$$h\in (\H_{g^{-1}(V)}, W) \ {\rm and} \ z\in g^{-1}(V).$$
Therefore, $z\in g^{-1}(V)\s h^{-1}(W
)$ and $e(h,z)=h(z)\in W$.

Thus, the evaluation map $e$ is continuous and, therefore, the
topology $t_{1,s}^Z$ is admissible. $\Box$

\begin{definition}
{\rm Let $Z$ be a space. A space $Y$ is called {\it
$Z$-corecompact}, if  for every $y\in Y$ and for every open
neighborhood $U$ of $y$ in $Y$, there exists a neighborhood
$f^{-1}(V)\in {\mathcal O}_Z(Y)$ of $y$ such that $f^{-1}(V)\s U$
and the subset $f^{-1}(V)$ is   $Z$-bounded in the space $U$ (in
symbols we write $f^{-1}(V)<<U$).}
\end{definition}

\begin{remark}
{\rm (1) In Definition 20 we considered the space $U$ to be a
subspace of the space $(Y,\tau_Y^Z)$, that is $U$ is the space which
is equipped with the  topology:
$$(\tau_Y^Z)_U=\{U\cap V:V\in \tau_Y^Z\}.$$ (2) Let $Y$ and $Z$ be two topological spaces, $U, W, U_i\in {\mathcal O}(Y)$, and  $f^{-1}(V)$, $f_i^{-1}(V_i)\in {\mathcal O}_Z(Y)$, where $i=1,\cdots,n$. We observe that:

\medskip

(i) If $f^{-1}(V)\s U<< W$, then $f^{-1}(V)<< W.$

\medskip

(ii) If $f^{-1}(V)<<  U\s W$, then $f^{-1}(V)<< W.$

\medskip

(iii) If $f_i^{-1}(V_i)<< U_i$, for every $i=1,\cdots,n$, then
$$\bigcup\{f_i^{-1}(V_i):i=1,\cdots,n\}<< \bigcup\{U_i:i=1,\cdots,n\}.$$}
\end{remark}

\begin{theorem}
{\rm Let $Z$ be a space and let $Y$ be a $Z$-corecompact space.
Then, the $t^Z_{1}$ topology on $C(Y,Z)$ is admissible.}
\end{theorem}

{\sl Proof.} It is sufficient to prove that the evaluation map:
$$e:C_{t^Z_{1}}(Y,Z)\times Y\to Z$$
is continuous.

For this, let $(f,y)\in C(Y,Z)\times Y$ and let $W \in {\mathcal
O}(Z)$ be such that:
$$e(f,y)=f(y)\in W.$$ Then, we have that:
$$y\in f^{-1}(W).$$ Since the space $Y$ is  $Z$-corecompact, there exists a neighborhood $g^{-1}(V)\in {\mathcal O}_Z(Y)$ of $y$, such
that:
$$g^{-1}(V)<< f^{-1}(W).$$ We now consider the set:
$$\H_{g^{-1}(V)}=\{U\in {\mathcal O}(Y): g^{-1}(V)<<U\}$$
and we prove that  the set $\H_{g^{-1}(V)} $ is open in the  $\tau_1^Z$ topology.

Indeed, let $f^{-1}(U)\in \H_{g^{-1}(V)}\cap {\mathcal O}_Z(Y)$,
$V_1\in {\mathcal O}(Y)$ and $f^{-1}(U)\s V_1$. Then, we have:
$$g^{-1}(V)<<f^{-1}(U)\s V_1$$ and, therefore, by Remark 2.5,  $g^{-1}(V)<< V_1$. Thus,  $V_1\in \H_{g^{-1}(V)}$. Now, let $\{f^{-1}_i(U_i):i\in I\}$ be a collection of sets of ${\mathcal O}_Z(Y)$, such
that:
$$\bigcup\{f^{-1}_i(U_i):i\in I\}\in \H_{g^{-1}(V)}$$ or
equivalently:
$$g^{-1}(V)<< \bigcup\{f^{-1}_i(U_i):i\in I\}.$$ Clearly, for every $y\in \bigcup\{f^{-1}_i(U_i):i\in
I\}$, there exists $i_y\in I$ such that $y\in
f^{-1}_{i_y}(U_{i_y})$. Since $Y$ is $Z$-corecompact, there exists a
neighborhood  $g_{i_y}^{-1}(V_{i_y})\in {\mathcal O}_Z(Y)$ of $y$,
such that:
$$y\in g^{-1}_{i_y}(V_{i_y})<< f^{-1}_{i_y}(U_{i_y}).$$ By all the
above we get:
$$ \bigcup\{g^{-1}_{i_y}(V_{i_y}):i_y\in I\}=\bigcup\{f^{-1}_{i_y}(U_{i_y}):i_y\in I\}
=\bigcup\{f^{-1}_i(U_{i}):i\in I\}.$$ Since $$g^{-1}(V)<<
\bigcup\{f^{-1}_i(U_{i}):i\in I\}$$ and since
$\{g^{-1}_{i_y}(V_{i_y}):i_y\in I\}$ is an open cover of
$\bigcup\{f^{-1}_i(U_i):i\in I\}$ with respect to the topology
$(\tau_Y^Z)_{\bigcup\{f^{-1}_i(U_i):i\in I\}}$, there exist finitely
many sets $g^{-1}_{i_{y_1}}(V_{i_{y_1}}),\cdots,
g^{-1}_{i_{y_n}}(V_{i_{y_n}})$ of this collection  such that:
$$g^{-1}(V)\s \bigcup\{g^{-1}_{{y_{i_k}}}(V_{{y_{i_k}}}):k=1,2,\cdots,n\}.\eqno(3)$$ We now consider the sets
$f^{-1}_{i_{y_1}}(U_{i_{y_1}}),\cdots, f^{-1}_{i_{y_n}}(U_{i_{y_n}})$
of ${\mathcal O}_Z(Y)$, for which:
$$g^{-1}_{i_{y_k}}(V_{i_{y_k}})<< f^{-1}_{i_{y_k}}(U_{i_{y_k}}), \ \ k=1,2,\cdots,n.$$ By Remark 2.5 and by (3)
above, we get:
 $$ g^{-1}(V)\s \bigcup\{g^{-1}_{i_{y_k}}(V_{i_{y_k}}):k=1,2,\cdots,n\}<<\bigcup\{f^{-1}_{i_{y_k}}(U_{i_{y_k}}):k=1,2,\cdots,n\}$$
and, therefore:
$$ g^{-1}(V)<<\bigcup\{f^{-1}_{i_{y_k}}(U_{i_{y_k}}):k=1,2,\cdots,n\}.$$ So,
$\bigcup\{f^{-1}_{i_{y_k}}(U_{i_{y_k}}):k=1,2,\cdots,n\}\in \H_{g^{-1}(V)}$.
Thus, the set $\H_{g^{-1}(V)}$ is open in the $\tau_1^Z$ topology.

Also, since $g^{-1}(V)<< f^{-1}(W)$, we have that $f\in
(\H_{g^{-1}(V)},W )$. Thus, the set $(H_{g^{-1}(V)}, W)\times
g^{-1}(V)$ is an open neighborhood of $(f,y)$ in
$C_{t_{1}^Z}(Y,Z)\times Y$.

We finally prove that:
$$e((\H_{g^{-1}(V)}, W)\times g^{-1}(V))\s W.$$ For this, let $(h,z)\in (\H_{g^{-1}(V)}, W
)\times g^{-1}(V)$.  Then, $h\in (\H_{g^{-1}(V)}, W
)$ and $z\in g^{-1}(V)$. Therefore, $z\in g^{-1}(V)<< h^{-1}(W
)$ and $e(h,z)=h(z)\in W$.

Thus, the evaluation map $e$ is continuous and, consequently, the
$t_1^Z$ topology is admissible. $\Box$

\medskip

By Theorem 2.9 and Remark 2.3 we obtain the following corollary.

\begin{corollary}
{\rm  Let $Z$ be any space and let    $Y$ be  a $Z$-corecompact
space. Then,
 the  $t_{1,s}^Z$ topology on $C(Y,Z)$
is admissible.}
\end{corollary}

\section{Admissible topologies on ${\mathcal O}_Z(Y)$}

Let $I\!\! H\subseteq {\cal O}_Z(Y)$, ${\cal H}\subseteq C(Y,Z)$
 and let
$U\in {\cal O}(Z)$. We set:
$$(I\!\! H, U)
=\{f\in C(Y,Z): f^{-1}(U)\in I\!\! H\}$$
and
$$({\cal H},U)=\{f^{-1}(U):f\in {\cal H}\}.$$

\begin{definition} {\rm (see \cite{GEO2})
(1) Let $\tau $ be a topology on ${\mathcal O}_Z(Y)$. The topology
on $C(Y,Z)$, for which the set

$$\{ (I\!\! H,U):I\!\! H\in\tau ,U\in{\cal O}(Z)\} $$

\noindent
is a subbasis,
is
 called the {\it
dual topology to $\tau $} and is denoted by $t(\tau )$.

\medskip
\noindent (2) Let $t$ be a topology on $C(Y,Z)$. The topology on
${\mathcal O}_Z(Y)$, for which the set:
$$
\{ ({\mathcal H},U):{\cal H}\in t, U\in{\mathcal O}(Z)\}$$ \noindent
is a subbasis, is called {\it the  dual topology to $t$} and is
denoted by $\tau (t)$.}
\end{definition}

Let $X$ be a space and let $G:X\to C(Y,Z)$ be a  map. By
$\overline{G}$ we denote the map from $X\times {\mathcal O}(Z)$ to
${\mathcal O}_Z(Y)$, for which $\overline{G}(x,U)=(G(x))^{-1}(U)$
for every $x\in X$ and $U\in {\mathcal O}(Z)$.

Let $\tau$ be a topology on  ${\mathcal O}_Z(Y)$. We say that a map
$M$  from $X\times {\mathcal O}(Z)$ to ${\mathcal O}_Z(Y)$ is {\it
continuous with respect to the
 first variable}  if for every fixed element
$U$ of ${\mathcal O}(Z)$, the map
$M_U:
X\to     ({\mathcal O}_Z(Y),\tau)$,
for  which $M_U(x)=M(x,U)$
for every $x\in X$,
 is continuous.

\begin{definition}{\rm (see \cite{GEO2})
A topology $\tau$ on ${\cal O}_Z(Y)$ is called {\it admissible}, if
for every space $X$ and for every map $G:X\to C(Y,Z)$ the continuity
 with respect to the first variable of the map
$\overline{G}:X\times {\mathcal O}(Z)\to ({\cal O}_Z(Y),\tau)$
implies the continuity of the map $\widetilde{G}:X\times Y\to Z$}
\end{definition}

It is known that (see \cite{GEO2}):

(1) A topology $t$ on $C(Y,Z)$ is admissible if and only if the
topology $\tau(t)$ on ${\mathcal O}_Z(Y)$ is admissible.

(2) A topology $\tau$ on ${\mathcal O}_Z(Y)$  is admissible if and only if the
topology $t(\tau)$ on  $C(Y,Z)$  is admissible.

\begin{theorem} {\rm The following statements hold:

(1) If $Y$ is a regular locally compact space (or a regular locally $Z$-compact space), then the topologies
$\tau (t^Z_{co})$, $\tau (t^Z_1)$, and $\tau (t^Z_{1,s})$ on
${\mathcal O}_Z(Y)$ are admissible.

(2) If $Y$ is a corecompact space (or a  $Z$-corecompact space), then the topologies
   $\tau(t^Z_1)$ and $\tau(t^Z_{1,s})$ on ${\mathcal O}_Z(Y)$  are admissible.

(3) If $Y$ is a locally bounded  space (or a locally $Z$-bounded  space), then the topology
 $\tau (t^Z_{1,s})$ on ${\mathcal O}_Z(Y)$ is admissible.}
 \end{theorem}

 {\sl Proof.} The proof of this theorem follows immediately from Theorems 2.5, 2.6, 2.7 and
 2.8. $\Box$

\begin{corollary} {\rm The following statements are true:

(1) If $Y$ is a regular locally compact space (or a regular locally
$Z$-compact space), then the topologies $t(\tau (t^Z_{co}))$,
$t(\tau (t^Z_1))$, and $t(\tau (t^Z_{1,s}))$ on $C(Y,Z)$ are
admissible.

(2) If $Y$ is a corecompact space (or a  $Z$-corecompact space), then the topologies
   $t(\tau(t^Z_1))$ and $t(\tau(t^Z_{1,s}))$ on $C(Y,Z)$  are admissible.

(3) If $Y$ is a locally bounded  space (or a locally $Z$-bounded  space), then the topology
 $t(\tau (t^Z_{1,s}))$ on $C(Y,Z)$ is admissible.}
 \end{corollary}

\begin{corollary} {\rm The following propositions are true:

(1) If $Y$ is a regular locally compact space (or a regular locally $Z$-compact space), then the topologies
$\tau(t(\tau (t^Z_{co})))$, $\tau(t(\tau (t^Z_1)))$, and $\tau(t(\tau (t^Z_{1,s})))$ on
${\mathcal O}_Z(Y)$ are admissible.

(2) If $Y$ is a corecompact space (or a  $Z$-corecompact space), then the topologies
   $\tau(t(\tau(t^Z_1)))$ and $\tau(t(\tau(t^Z_{1,s})))$ on ${\mathcal O}_Z(Y)$  are admissible.

(3) If $Y$ is a locally bounded  space (or a locally $Z$-bounded  space), then the topology
 $\tau(t(\tau (t^Z_{1,s})))$ on ${\mathcal O}_Z(Y)$ is admissible.}

 \end{corollary}

\section{Some open questions}

In this section we give some interesting in our opinion open
questions applied to the topologies $t_{co}^Z$,
 $t_1^Z$, and $t_{1,s}^Z$.

\medskip

{\bf Question 1.} Let $Y$ and $Z$ be two topological spaces. Is the topology $t_1^Z$ on $C(Y,Z)$ regular  in the case where
$Z$ is regular?

\medskip

{\bf Question 2.} Let $Y$ and $Z$ be two topological spaces. Is the topology $t_1^Z$ on $C(Y,Z)$  completely regular in the case where $Z$ is  completely regular?

\medskip

{\bf Question 3.} Find two topological spaces $Y$ and $Z$ such that

{\bf 3.1.} $t_{co}^Z\not =t_{co}$.

{\bf 3.2.} $t_1^Z\not =t_{Is}$.

{\bf 3.3.} $t_{1,s}^Z\not =t_{sIs}$.

\medskip

{\bf Question 4.}
 Do the topologies $t_{co}$ and $t_{co}^Z$ coincide on the set
$C(\mathbb R^{\omega},\mathbb R)$ where $\mathbb R$ is  the set of real numbers  with the usual topology
and $\omega$ is the first infinite cardinal?

\medskip

{\bf Question 5.}
 Do the topologies $t_{co}$ and $t_{co}^Z$ coincide on the set
 $C(\mathbb N^{\omega},\mathbb N)$, where  $\mathbb N$ is the set of natural numbers  with its usual topology?\

\medskip

{\bf Question 6.} Let $X$, $Y$, and $Z$ be three topological spaces.
Is the map: $$T: C_{t^Y_{1,s}}(X,Y) \times C_{t^Z_{1,s}}(Y,Z) \to
C_{t^Z_{1,s}}(X,Z),$$ with $T(f,g)=g\circ f$, for every $(f,g) \in
C(X,Y)\times C(Y,Z)$, continuous in the case where $Y$ is locally
$Z$-bounded?

\medskip

{\bf Question 7.} Let $X$, $Y$, and $Z$ be three topological spaces.
Is the map: $$T: C_{t^Y_{1}}(X,Y) \times C_{t^Z_{1}}(Y,Z) \to
C_{t^Z_{1}}(X,Z),$$ with $T(f,g)=g\circ f$, for every $(f,g) \in
C(X,Y)\times C(Y,Z)$, continuous in the case where $Y$ is
$Z$-corecompact?

\medskip

{\bf Notation.} Let $X$ be a space and $F:X\times Y\to Z$ be a
continuous map. By $F_x$, where $x\in X$, we denote the continuous
map of $Y$ into $Z$ such that $F_x(y)=F(x,y)$, $y\in Y$. By
$\widehat{F}$ we denote the map of $X$ into the  set $C(Y,Z)$ such
that $\widehat{F}(x)=F_x$, $x\in X$.

We recall that a topology $t$ on $C(Y,Z)$ is called {\it splitting}
 if for every space $X$, the continuity of a map
$F:X\times Y\to Z$  implies that of the
 map
$\widehat{F}:X\to C_{t}(Y,Z)$ (see \cite{ARE} and \cite{DUG}).

It is known that:

(1)  The compact open  topology $t_{co}$ is always splitting (see
\cite{ARE} and \cite{FOX}).

(2) The Isbell topology is always splitting (see \cite{ISBE},
\cite{LAW1}, \cite{MN}, and \cite{SW}).

(3)  If $Z$ is the Sierpi\'{n}ski space   and $Y$ is an arbitrary space,
then  the Isbell
 topology coincides with the greatest splitting topology  (see \cite{PAP} and
 \cite{SW}).

By Remark 2.3(2) if $Z$ is the Sierpi\'{n}ski space, then the topology
$t^Z_1$ coincides with the Isbell topology on $C(Y,Z)$ and,
therefore, this topology is a splitting topology. In this case the
topology $t^Z_{co}$ is also splitting.

\medskip

{\bf Question 8.} Find necessary and sufficient conditions for the
space $Z$ such that the topology $t^Z_1$ on $C(Y,Z)$ to be
splitting.

\medskip

{\bf Question 9.} Find necessary and sufficient conditions for the
space $Z$ such that the topology $t^Z_{co}$ on $C(Y,Z)$ to be
splitting.

\medskip

{\bf Question 10.} Let $Y$ be an arbitrary topological space and
$Z=\{0,1\}$ with the discrete topology.  Which of the
following relations is true or false on $C(Y,Z)$?

(i) $t_{co}\not = t^Z_{co}$,

(ii) $t_{Is}\not = t^Z_1$, and

(iii) $t_{sIs} \not = t_{1,s}^Z$.

\medskip

{\bf Question 11.} Let $Y$ be an arbitrary topological space and
$Z=\{0,1, {1\over 2}, {1\over 3}, {1\over 4},\ldots\}$  with the usual
topology of the real line. Which of the following relations is true or false on $C(Y,Z)$?

(i) $t_{co}\not = t^Z_{co}$,

(ii)  $t_{Is}\not = t^Z_1$, and

(iii) $t_{sIs} \not = t_{1,s}^Z$.

\medskip

Let $Y$ be a corecompact space which is not basic locally compact.
(A space is basic locally compact if for every point there exists a basis of
compact neighbourhoods). Then, the Isbell topology $t_{Is}\equiv
t_1^{\bf S}$, in $C(Y,{\bf S})$ does not coincide to the compact
open topology $t_{co}\equiv t_{co}^{\bf S}$ (see \cite{PAP} and
\cite{SW}).

\medskip

{\bf Question 12.} For what spaces $Y$ and $Z$ does the equality
$t_1^Z=t_{co}^Z$ hold?

\bigskip
\noindent
{\bf Acknowledgements.} The authors would like to thank the referee for very helpful comments and suggestions.

\bigskip

$\begin{array}{ll}
\rm{Email \ addresses:} & {\rm georgiou@math.upatras.gr \ (Dimitris \ Georgiou)}\\
 & {\rm megariti@master.math.upatras.gr \ (Athanasios \ Megaritis)}\\
 & {\rm KXP878@bham.ac.uk \ (Kyriakos \ Papadopoulos)}
\end{array}$

\end{document}